\documentclass[final,5p,times,twocolumn]{elsarticle}

\usepackage{amssymb}
\usepackage{float}
\usepackage{amsmath}
\usepackage{url}
\usepackage{hyperref}
\usepackage{indentfirst}
\usepackage{cases}
\usepackage{epstopdf}
\usepackage{multirow}
\usepackage{amsfonts}
\usepackage{amssymb}
\usepackage{extarrows}
\usepackage{graphicx}
\usepackage{amsthm}

\newtheorem{lemma}{Lemma}
\newtheorem{definition}{Definition}
\hypersetup{colorlinks = true,linkcolor = blue,anchorcolor =red,citecolor = blue,filecolor = red,urlcolor = red,
            pdfauthor=author}

\journal{Fundamental Research, 4 (2024) 841-844, by L.~Wang \& B.-Q.~Ma, \href{https://doi.org/10.1016/j.fmre.2023.01.002}{doi:10.1016/j.fmre.2023.01.002}}

\begin{document}

\begin{frontmatter}



\title{A concise proof of Benford's law}


\author[inst1]{Luohan Wang}

\affiliation[inst1]{organization={School of Physics},
            addressline={Peking University}, 
            city={Beijing},
            postcode={100871}, 
            country={China}}

\author[inst1,inst2,inst3]{Bo-Qiang Ma\corref{cor1}}
\cortext[cor1]{mabq@pku.edu.cn}
\affiliation[inst2]{organization={Center for High Energy Physics},
            addressline={Peking University}, 
            city={Beijing},
            postcode={100871}, 
            country={China}}

\affiliation[inst3]{organization={Collaborative Innovation Center of Quantum Matter},
            city={Beijing},
            country={China}}

\begin{abstract}
This article presents a concise proof of the famous Benford's law when the distribution has a Riemann integrable probability density function and provides a criterion to judge whether a distribution obeys the law. The proof is intuitive and elegant, accessible to anyone with basic knowledge of calculus, revealing that the law originates from the basic property of the human number system.
The criterion can bring great convenience to the field of fraud detection. 
\end{abstract}





\end{frontmatter}


\section{Introduction}

One may instinctively assume that the occurrence of the first digit of a randomly selected number is uniformly distributed among 1 to 9.  Nonetheless, the Newcomb-Benford law, also known as the first-digit law or Benford's law, reveals a completely different fact. In many natural and human phenomena, the distribution of the first significant digit of a random number follows a logarithm-type law.

Early in 1881, 
Newcomb \cite{Newcomb1881} noticed that the first pages of the logarithmic table tend to wear out faster than the last ones. Consequently, he found that the first significant figure is oftener 1 than any other digit, with the frequency diminishing from 1 to 9, and figured out the logarithm distribution of the first two significant digits of natural numbers. In 1938, by analyzing over 20,000 numbers collected from 20 widely divergent sources, 
Benford \cite{Benford1938} found that these data show unbelievable adherence to a logarithm distribution, and rediscovered the empirical law
\begin{equation}
    P(d)=\log_{10}\left(1+\frac{1}{d}
    \right),\quad d=1,2,\cdots,9, \label{eq1}
\end{equation}
where $P(d)$ is the probability of a randomly selected number with a first digit of $d$. However, no convincing proof was proposed by him at that time.

This bias in the frequency of occurrence of the leading digits of numbers is evident in everything from the lengths of rivers, the areas of countries, stock prices, and population numbers to death rates. In terms of mathematics, Benford's law appears in integer sequences, such as the celebrated Fibonacci sequence \cite{TPHill2015} and the factorial sequence \cite{MillerSJ2015}, and the distribution of prime numbers follows a general Benford's law \cite{luque2009}. The law also emerges in the field of physics. 
  The collection of many physics constants \cite{Burke1991} conforms to this law. The half-lives of unstable nuclei \cite{NDD2008,Ni2009,LXJ2011} and many other experimental quantities \cite{Jiang2011} in nuclear physics also agree with the logarithm distribution. Its applications in physics involve astronomy \cite{SHAO2010feb}, dynamic system \cite{Snyder2001}, statistical physics \cite{SHAO2010apr,Shao2012oct}, nuclear physics~\cite{NDD2008,Ni2009,LXJ2011,Jiang2011} and particle physics~\cite{Shao2009}. For example, the Bose-Einstein distribution conforms to it strictly whatever the temperature is, while the Fermi-Dirac distribution fluctuates slightly in a periodic manner depending on the temperature of the system \cite{SHAO2010apr}.

The detection of frauds is one of the most prominent applications of the Newcomb–Benford law for significant digits. In 1972, a seminal letter by 
Varian \cite{Letters} suggested the notion that Benford's law can be applied in detecting fraud data. In 1992, 
Nigrini proposed in his doctoral dissertation to use the law to check for false accounting. Since then, the law has become rather popular in auditing financial data, testing the fairness of electoral processes \cite{Luis2011}, while a study anticipating the validity of this law in detecting frauds in international trades was carried out in 2019 \cite{Parma2019}. The law also succeeds in application to many other domains, including social science \cite{Leemis2000} and environmental science \cite{Richard2005}.

Nevertheless, Benford's law has its scope of application. It works in a great many natural date sets and some human phenomena, but most artificial data sets and distributions that do not span several orders of magnitude (telephone number, height, weight) or are artificially truncated do not follow Benford's law. 

With regard to the underlying mechanism of Benford's law, there have been explanations from different views. However, due to the intriguing feature of Benford's law, 
the available explanations are not completely satisfactory. 
They are either under a restricted condition or mathematically not rigorous. The inchoate research focused on the cases where the samples are number sets defined on $\mathbb{R^+}$, such as the set of natural integers. In 1976, 
Cohen \cite{Cohen1976} raised the condition under which a measure on $\mathbb{N}$ satisfies Benford's law. In 1992, 
Jech \cite{JECH1992} stated the necessary and sufficient conditions for a probability measure with a finitely additive function defined on $\mathbb{R^+}$ to observe the law. The law also has proof of physical versions due to its frequent appearance in physics. 
Lemon \cite{Lemons2019} and 
Iafrate \cite{Iafrate2015} derived this law using the maximum entropy model in thermodynamics.

An explanation widely considered as rigorous is by 
Hill \cite{hill1995} in 1995. 
He illustrated that when the distributions and the sample data are both selected at random, the first digits of the combinations follow Benford's law. However, this work did not answer the question of why data from a specific distribution obey this law and did not pay much attention to the internal reasons and conditions of the law.

In addition, little attention has been paid to the circumstances where the random variable set is continuous and has a continuous probability density function defined on it. Studying this topic helps to explain why a continuous distribution follows or violates the law. For example, among the common distributions, the exponential distribution conforms well to the law, while the uniform distribution and the normal distribution 
do not obey the law. Research conducted by 
Engel and 
Leuenberger \cite{ENGEL2003} proved that exponential functions obey this law within an allowed error range. In 2019, a study at Peking University gave a proof of the law when the probability density function has an inverse Laplace transform \cite{proof,proof2}, revealing that the first digit law originates from the basic property of the human number system.

In this article, first, we aim to derive this law for distributions with Riemann integrable probability density functions on the positive real number set, and the basic idea of the proof is inspired by Refs.~\cite{proof,proof2}. Our proof is concise, intuitive, and accessible to anyone with basic knowledge of calculus. Second, we estimate the potential error caused by a hypothesis in the proof and discuss the criterion to judge the validity of this law in a distribution.

\section{Proof for the continuous case}

\subsection{The general proof of Benford's law}

Let the sample space of the variable $x$ be the positive real number set $\mathbb{R^+}$, and let $f(x)$ be a continuous probability density function of an arbitrary distribution defined on $\mathbb{R^+}$, with $f(x)$ satisfying the normalization condition. The set of positive numbers with the first significant of $d\ (d=1,2,\cdots,9)$ is
\begin{equation}
S(d)=\bigcup\limits_{n=-\infty}^{+\infty}[d\cdot 10^n,(d+1)\cdot 10^n).
\end{equation}
Therefore, it is obvious that
\begin{equation}\label{eq(Pd)}
P(d)=\sum_{n=-\infty}^{+\infty}\int_{d\cdot 10^n}^{(d+1)\cdot 10^n}f(x)\,\text{d}x,
\end{equation}
where $P(d)$ is the probability of a number (base 10) chosen randomly from the distribution $f(x)$ beginning with $d$. Evidently, $P(d)$ can also be written as
\begin{equation}
    P(d)=\sum_{n=-\infty}^{+\infty}\int_{d\cdot 10^n}^{(d+1)\cdot 10^n}f(x)\,\text{d}x\,\Delta n ,\quad (\Delta n=1).
\end{equation}

Then we will make an important assumption that the sum in the above equation is approximately equal to the integral as follows, which is similar to Eq.~(11) in Refs.~\cite{proof,proof2}
\begin{align}\label{approx}
P(d)\approx \int_{-\infty}^{+\infty}\left[\int_{d\cdot10^n}^{(d+1)\cdot10^n}f(x)\,\text{d}x\right]\text{d}n,
\end{align}
since $\Delta n$ is not exactly the differential d$n$, this equation does not always hold. This explains why some $f(x)$ break Benford's law when the difference between the two equations cannot be ignored. We discuss this situation in the next section.

Substitute the variable $n$ with $t=d\cdot 10^n$, $x$ with $x=ty$, therefore $\text{d}n=\frac{1}{t\ln 10}\,\text{d}t$, and the same $P(d)$ in Eq.~(\ref{approx}) can be written as
\begin{align}
P(d)=\frac{1}{\ln 10}\int_{0}^{+\infty}\text{d}t\int_{1}^{1+\frac{1}{d}}f(ty)\,\text{d}y.
\end{align}
\begin{definition}
$\phi(y,t)=f(ty)$, and $\phi(y,t)$ is continuous on\\ $\left[1,1+\frac{1}{d}\right]\times[0,+\infty)$.
\end{definition}
Apparently, the integral $\int_{0}^{+\infty}\phi(y,t)\,\text{d}t$ is uniform convergent on $\left[1,1+\frac{1}{d}\right]$.
\begin{lemma}\label{lemma}
If $\phi(y,t)$ is continuous on $[a,b]\times[c,+\infty)$, $a,b,c \in \mathbb{R^+}$, and the integral $\int_{c}^{+\infty}\phi(y,t)\,\mathrm{d} t$ is uniform convergent on $[a,b]$, then there is
\begin{equation}   \int_{c}^{+\infty}\mathrm{d}t\int_{a}^{b}\phi(y,t)\,\mathrm{d}y=\int_{a}^{b}\mathrm{d}y\int_{c}^{+\infty}\phi(y,t)\,\mathrm{d}t,
\end{equation}
which is a basic theorem in calculus.
\end{lemma}

Therefore, according to Lemma~\ref{lemma},
\begin{align}
    P(d)&=\frac{1}{\ln 10}\int_{0}^{+\infty}\text{d}t\int_{1}^{1+\frac{1}{d}}\phi(y,t)\,\text{d}y \notag\\
&=\frac{1}{\ln 10}\int_{1}^{1+\frac{1}{d}}\text{d}y\int_{0}^{+\infty}\phi(y,t)\,\text{d}t \notag\\
&=\frac{1}{\ln 10}\int_{1}^{1+\frac{1}{d}}\frac{1}{y}\,\text{d}y\int_{0}^{+\infty}f(x)\,\text{d}x \notag\\
&=\log_{10}\left(1+\frac{1}{d}\right) \label{pfofNBL}
\end{align}
is exactly the logarithm relation in Benford's law. The proof above reveals that Benford's law is essentially a mathematical principle. It arises simply because of humans' method of counting.

\subsection{The extensional versions of Benford's law}

Apart from the commonest first digit law in the decimal system, Benford's law has other properties and extensional versions. The law is both scale-invariant and base-invariant. The scale-invariance is the direct conclusion of the proof above. By replacing $x$ with $ax$ ($a$ is a positive real constant) in Eq.~(\ref{eq(Pd)}), we can find that the result is the same since the integrating range is from zero to infinity.

In 1995, 
Hill \cite{hill1995} derived a general $i$th-significant digit law: For all positive integers $k$, all $d_1\in\left\{1,2,\cdots,9\right\}$ and all $d_j\in\left\{0,1,\cdots,9\right\}$, $j=1,2,\cdots,k$, the probability of a number beginning with first $k$ digit, $D_j=d_j$ can be expressed as
\begin{align}\label{generalNBL}
    P&\left(D_1=d_1,\cdots,D_k=d_k\right)\notag\\
    &=\log_{10}\left[1+\left(\sum_{i=1}^{k}d_i\cdot 10^{k-i}\right)^{-1}\right].
\end{align}

Here is an extension of Eq.~(\ref{generalNBL}) to a number with continuous coefficient in base $m$.
Suppose $x=k\cdot m^n,\ m\in \mathbb{N^+},\ k \in [1,m)$, and $n$ is an integer. Specially, this is the scientific notation of $x$ when $m=10$. Then use the same method mentioned above, we can calculate that the probability of a number $x$ with its $k\in[a,b),\ a\geq1,\ b<m$ is
\begin{align}
    P_m(a,b)&=\sum_{n=-\infty}^{+\infty}\int_{a\cdot m^n}^{b\cdot m^n}f(x)\,\text{d}x \notag \\
    &=\log_m b/a. \label{generalNBL2}
\end{align}
 Clearly, this general version of Benford's law is base-invariant and Eq.~(\ref{generalNBL}) is a special case when $m=10,\ a=\sum\limits_{i=1}\limits^{k}d_i\cdot 10^{k-i},\ b=a+1$ of Eq.~(\ref{generalNBL2}).

\section{Discussions over the errors}

In the proof above, we use an approximation as Eq.~(\ref{approx}). Apparently, replacing the sum with an integration is not universally feasible. When $n$ gets large, the integral interval of $f(x)$ for each $n$ grows wider, and the integration is more possibly to disagree with the original sum. Therefore, we should discuss carefully the rationality of this replacement, as well as  what kind of $f(x)$ obeys Benford's law. 
\begin{definition}
$g_d(x)$ is a function satisfying
\begin{flalign}
    g_d(x)=\sum_{n=-\infty}^{+\infty}\left[\eta\left(x-d\cdot 10^n\right)-\eta\left(x-(d+1)\cdot 10^n\right)\right],
\end{flalign}
which is the same function in Refs.~\cite{proof,proof2},
and $\eta(x)$ is the Heaviside step function
\begin{align}
\eta(x)=\left\{
\begin{aligned}
0\quad x<0 \\
1\quad x\geq 0
\end{aligned}
\right.
.
\end{align}
\end{definition}
Therefore, $P(d)$ can also be rewritten as
\begin{equation}
    P(d)=\int_0^{+\infty}g_d(x)f(x)\,\text{d}x.
\end{equation}
\begin{definition}
Define a function $\Delta(x)$ as follows:
\begin{equation}
    \Delta(x)=g_d(x)-\log_{10}\left(1+\frac{1}{d}\right).
\end{equation}
Denote the Error function $\mathrm{Er}(f)$ a functional of $f(x)$ which represents the difference between the original sum and the integration in Eq.~(\ref{approx}), therefore
\end{definition}
\begin{align}\label{Er}
    \mathrm{Er}(f)&=\int_0^{+\infty}g_d(x)f(x)\,\text{d}x-\log_{10}\left(1+\frac{1}{d}\right)\notag
    \\&=\int_0^{+\infty}f(x)\left[g_d(x)-\log_{10}\left(1+\frac{1}{d}\right)\right]\text{d}x\notag
    \\&=\int_0^{+\infty}f(x)\Delta(x)\,\text{d}x,
\end{align}
$\Delta(x)$ oscillates within $\left[-\log_{10}\left(1+\frac{1}{d}\right),1-\log_{10}\left(1+\frac{1}{d}\right)\right]$ and the interval is proportional to $10^n,\ n\in \mathbb{Z}$. Qualitatively, if $f(x)$ strides across many orders of magnitude, the positive and negative part of $\Delta(x)$ may be offset. However, if $f(x)$ is defined on a narrow interval, $\mathrm{Er}(f)$ could be significant, that's why Benford's law is broken when the distribution is narrow. It can also explain why artificial data sets always contradict this law, since humans tend to make a uniform distribution while uniform distributions violate Benford's law. 

To judge whether the replacement used in Eq.~(\ref{approx}) is allowed, we only have to calculate the integration in Eq.~(\ref{Er}). If $\mathrm{Er}(f)$ is minute, $f(x)$ obeys Benford's law. The examples in Refs.~\cite{proof,proof2} and Ref.~\cite{ENGEL2003} both verify the effectiveness of this criterion. Here we take a simple case where $f(x)=\lambda e^{-\lambda x},\ \lambda>0$ for example:
\begin{equation}
\text{Er}(f)=\sum\limits_{n=-\infty}\limits^{+\infty}e^{-\lambda d \cdot10^n}\left(1-e^{-\lambda\cdot 10^n}\right)-\log_{10}\left(1+\frac{1}{d}\right).
\end{equation}
According to the calculation in Ref.~\cite{ENGEL2003}, $|\text{Er}(f)|\leq 0.03$ for any $\lambda$ and $d$. To give an intuitive illustration, we plot a figure of the difference between Benford's distribution and the exact distributions of exponential functions as shown in Fig.1. 
\begin{figure}[H]
    \centering
    \label{Fig1}
    \includegraphics[width=0.5\textwidth]{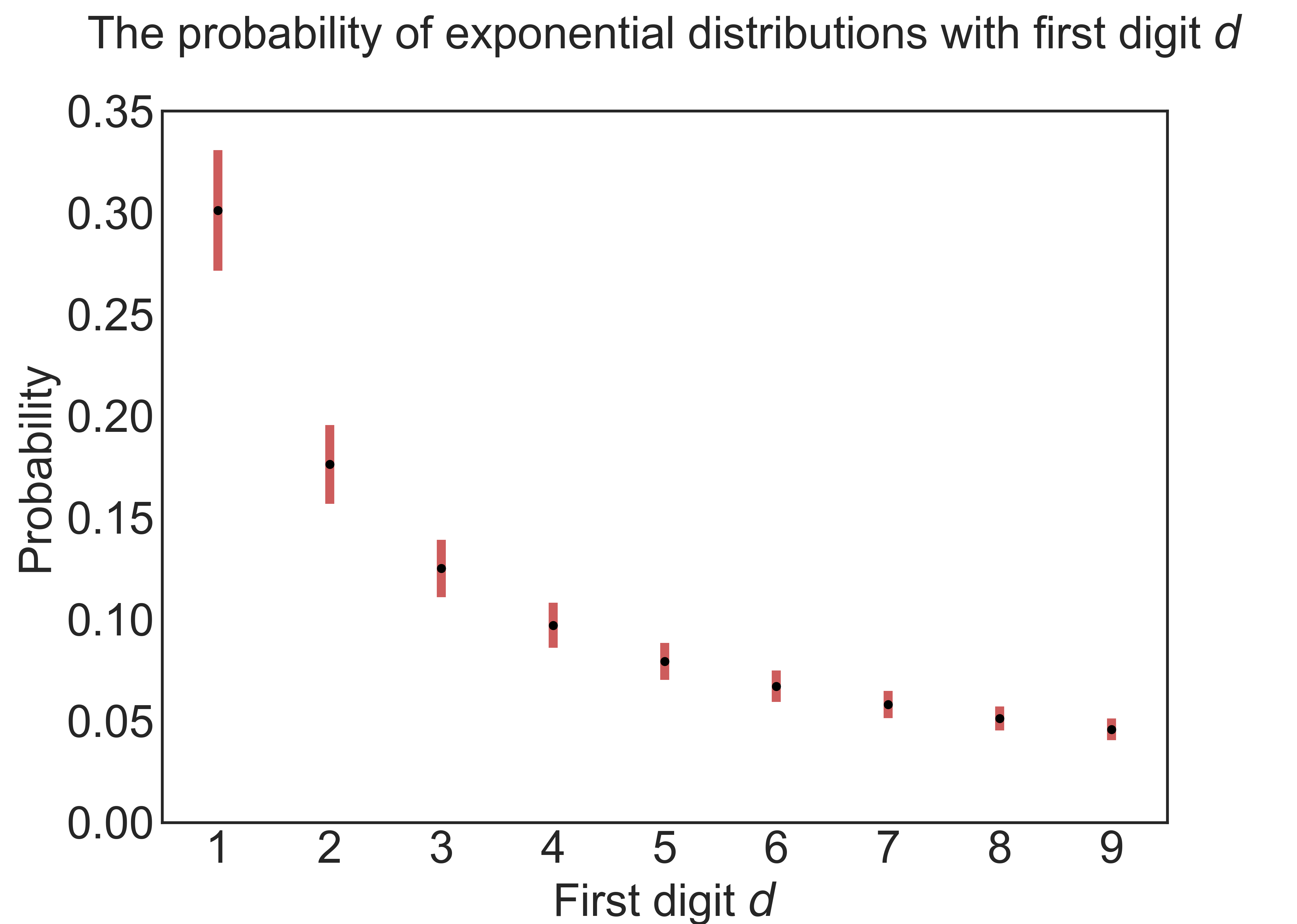}
    \caption{Graph of $P(d)$ and the range of error.
    The red part is the maximum $\vert \text{Er}(f)\vert$ of $f(x)$ with different $\lambda$ for each $d$, and the black dots are Benford's distributions.}
\end{figure}

Apparently, Benford's law is a good estimate for exponential distributions. This also explains the prevalence of Benford's law in nature. A number of distributions in nature are exponential, or can be written as superposition of exponential functions, such as $f(x)=\int_0^{+\infty} c(t)e^{-tx}\text{d}t$. The normalized probability density function $f(x)$ can be seen as the superposition of $e^{-tx}$ with coefficient $c(t)$ in this expression, meaning $f(x)$ has an inverse Laplace transformation. Therefore, it is intuitive to guess that such $f(x)$ with positive $c(t)$ converges to Benford's law, which is exactly the conclusion of Refs.~\cite{proof,proof2}.

What else, though the proof above is under the condition that $f(x)$ is continuous, it can be extended to discontinuous $f(x)$. If $f(x)$ can be approximated by a continuous function, in other words, for all ${\displaystyle \varepsilon > 0}$, there exists a continuous function $h(x)$ on $[0,+\infty]$ satisfying
\begin{equation}
    \int_0^{+\infty}\left|f(x)-h(x)\right|\text{d}x < \varepsilon,
\end{equation}
then the difference between $\mathrm{Er}(f)$ and $\mathrm{Er}(h)$ is also infinitesimal. As long as $\mathrm{Er}(h)$ is a small number, $f(x)$ converges to the logarithm distribution. $h(x)$ can be the Fourier expansion of $f(x)$ or a  piecewise linear function, which is easy to calculate its $\mathrm{Er}(h)$. Therefore, the condition of the proof is weakened as $f(x)$ is Riemann integrable on $\mathbb{R^+}$.

We also have to mention that our methodology does not work for all the distributions. The method loses efficacy when $f(x)$ is not Riemann integrable or the distribution does not have an asymptotic density. Related works in these circumstances can be found in Refs.~\cite{Cohen1976,Cigler1964,Flehinger1966}.

\section{Summary}

Over more than one hundred years, Benford's law has drawn much attention in different disciplines. Despite its long history, the law still lacks a satisfactory explanation and is even considered as a mystery of nature. Does nature favor logarithms? Does nature count at the base of natural exponential functions instead of our natural number scale? Benford's law has inspired so many conjectures. 

In this article, we prove Benford's law for Riemann integrable probability density functions defined on $\mathbb{R^+}$ and give the criterion to estimate the accuracy of the law.
Our proof is concise, intuitive, and accessible to anyone with basic knowledge of calculus.
With the criterion we can easily determine whether a distribution with a continuous variable obeys Benford's law, thus such criterion can bring great convenience to the field of fraud detection. 



\section*{Declaration of competing interest}
The authors declare that they have no conflict of interests.

\section*{Acknowledgments}
This work is supported by National Natural Science Foundation of China (Grant No.~12075003).

 \bibliographystyle{elsarticle-num} 
 \bibliography{ref.bib}

\begin{thebibliography}{10}
\expandafter\ifx\csname url\endcsname\relax
  \def\url#1{\texttt{#1}}\fi
\expandafter\ifx\csname urlprefix\endcsname\relax\def\urlprefix{URL }\fi
\expandafter\ifx\csname href\endcsname\relax
  \def\href#1#2{#2} \def\path#1{#1}\fi

\bibitem{Newcomb1881}
S.~Newcomb, Note on the frequency of use of the different digits in natural
  numbers, Amer. J. Math. 4~(1-4) (1881) 39--40.
\newblock \href {https://doi.org/10.2307/2369148} {\path{doi:10.2307/2369148}}.

\bibitem{Benford1938}
F.~Benford, The law of anomalous numbers, Proceedings of the American
  Philosophical Society 78~(4) (1938) 551--572,
  \url{http://www.jstor.org/stable/984802}.

\bibitem{TPHill2015}
A.~Berger, T.~P. Hill, An introduction to {Benford's} law, Princeton University
  Press, 2015, \url{https://www.jstor.org/stable/j.ctt1dr35m0 }.

\bibitem{MillerSJ2015}
P.~Diaconis, {The distribution of leading digits and uniform distribution mod
  1}, The Annals of Probability 5~(1) (1977) 72 -- 81.
\newblock \href {https://doi.org/10.1214/aop/1176995891}
  {\path{doi:10.1214/aop/1176995891}}.

\bibitem{luque2009}
B.~Luque, L.~Lacasa, The first-digit frequencies of prime numbers and riemann
  zeta zeros, Proceedings of the Royal Society A: Mathematical, Physical and
  Engineering Sciences 465~(2107) (2009) 2197--2216.
\newblock \href {https://doi.org/10.1098/rspa.2009.0126}
  {\path{doi:10.1098/rspa.2009.0126}}.

\bibitem{Burke1991}
J.~Burke, E.~Kincanon, Benford’s law and physical constants: The distribution
  of initial digits, American Journal of Physics 59~(10) (1991) 952--952.
\newblock \href {https://doi.org/10.1119/1.16838} {\path{doi:10.1119/1.16838}}.

\bibitem{NDD2008}
D.~Ni, Z.~Ren, {Benford’s} law and half-lives of unstable nuclei, The
  European Physical Journal A 38 (2008) 251--255.
\newblock \href {https://doi.org/10.1140/epja/i2008-10680-8}
  {\path{doi:10.1140/epja/i2008-10680-8}}.

\bibitem{Ni2009}
D.~Ni, L.~Wei, Z.~Ren, {Benford's} law and $\beta$-decay half-lives,
  Communications in Theoretical Physics 51~(4) (2009) 713--716.
\newblock \href {https://doi.org/10.1088/0253-6102/51/4/25}
  {\path{doi:10.1088/0253-6102/51/4/25}}.

\bibitem{LXJ2011}
X.~J. Liu, X.~P. Zhang, D.~D. Ni, Z.~Z. Ren, {Benford’s} law and
  cross-sections of {A}(n,$\alpha$){B} reactions, The European Physical Journal
  A 47 (2011) 78.
\newblock \href {https://doi.org/10.1140/epja/i2011-11078-3}
  {\path{doi:10.1140/epja/i2011-11078-3}}.

\bibitem{Jiang2011}
H.~Jiang, J.~Shen, Y.~Zhao, Benford{\textquotesingle}s law in nuclear structure
  physics, Chinese Physics Letters 28~(3) (2011) 032101.
\newblock \href {https://doi.org/10.1088/0256-307x/28/3/032101}
  {\path{doi:10.1088/0256-307x/28/3/032101}}.

\bibitem{SHAO2010feb}
L.~Shao, B.-Q. Ma, Empirical mantissa distributions of pulsars, Astroparticle
  Physics 33~(4) (2010) 255--262.
\newblock \href {https://doi.org/10.1016/j.astropartphys.2010.02.003}
  {\path{doi:10.1016/j.astropartphys.2010.02.003}}.

\bibitem{Snyder2001}
M.~A. Snyder, J.~H. Curry, A.~M. Dougherty, Stochastic aspects of
  one-dimensional discrete dynamical systems: Benford's law, Phys. Rev. E 64
  (2001) 026222.
\newblock \href {https://doi.org/10.1103/PhysRevE.64.026222}
  {\path{doi:10.1103/PhysRevE.64.026222}}.

\bibitem{SHAO2010apr}
L.~Shao, B.-Q. Ma, The significant digit law in statistical physics, Physica A:
  Statistical Mechanics and its Applications 389~(16) (2010) 3109--3116.
\newblock \href {https://doi.org/10.1016/j.physa.2010.04.021}
  {\path{doi:10.1016/j.physa.2010.04.021}}.

\bibitem{Shao2012oct}
L.~Shao, B.-Q. Ma, First-digit law in nonextensive statistics, Phys. Rev. E 82
  (2010) 041110.
\newblock \href {https://doi.org/10.1103/PhysRevE.82.041110}
  {\path{doi:10.1103/PhysRevE.82.041110}}.

\bibitem{Shao2009}
L.~Shao, B.-Q. Ma, First digit distribution of hadron full width, Modern
  Physics Letters A 24~(40) (2009) 3275--3282.
\newblock \href {https://doi.org/10.1142/S0217732309031223}
  {\path{doi:10.1142/S0217732309031223}}.

\bibitem{Letters}
HR.Varian, Letters to the editor, The American Statistician 26~(3) (1972)
  62--65.
\newblock \href {https://doi.org/10.1080/00031305.1972.10478934}
  {\path{doi:10.1080/00031305.1972.10478934}}.

\bibitem{Luis2011}
L.~Pericchi, D.~Torres, {Quick anomaly detection by the {Newcomb–Benford}
  law, with applications to electoral processes data from the USA, Puerto Rico
  and Venezuela}, Statistical Science 26~(4) (2011) 502 -- 516.
\newblock \href {https://doi.org/10.1214/09-STS296}
  {\path{doi:10.1214/09-STS296}}.

\bibitem{Parma2019}
A.~Cerioli, L.~Barabesi, A.~Cerasa, M.~Menegatti, D.~Perrotta,
  Newcomb--{Benford} law and the detection of frauds in international trade,
  Proceedings of the National Academy of Sciences 116~(1) (2019) 106--115.
\newblock \href {https://doi.org/10.1073/pnas.1806617115}
  {\path{doi:10.1073/pnas.1806617115}}.

\bibitem{Leemis2000}
L.~M. Leemis, B.~W. Schmeiser, D.~L. Evans, Survival distributions satisfying
  {Benford's} law, American Statistician 54~(4) (2000) 236--241.
\newblock \href {https://doi.org/10.1080/00031305.2000.10474554}
  {\path{doi:10.1080/00031305.2000.10474554}}.

\bibitem{Richard2005}
R.~J.~C. Brown, Benford{'}s law and the screening of analytical data: the case
  of pollutant concentrations in ambient air, Analyst 130 (2005) 1280--1285.
\newblock \href {https://doi.org/10.1039/B504462F}
  {\path{doi:10.1039/B504462F}}.

\bibitem{Cohen1976}
D.~I. Cohen, An explanation of the first digit phenomenon, Journal of
  Combinatorial Theory, Series A 20~(3) (1976) 367--370.
\newblock \href {https://doi.org/10.1016/0097-3165(76)90032-7}
  {\path{doi:10.1016/0097-3165(76)90032-7}}.

\bibitem{JECH1992}
T.~Jech, The logarithmic distribution of leading digits and finitely additive
  measures, Discrete Mathematics 108~(1) (1992) 53--57.
\newblock \href {https://doi.org/10.1016/0012-365X(92)90659-4}
  {\path{doi:10.1016/0012-365X(92)90659-4}}.

\bibitem{Lemons2019}
D.~S. Lemons, Thermodynamics of {Benford's} first digit law, American Journal
  of Physics 87~(10) (2019) 787--790.
\newblock \href {https://doi.org/10.1119/1.5116005}
  {\path{doi:10.1119/1.5116005}}.

\bibitem{Iafrate2015}
J.~R. Iafrate, S.~J. Miller, F.~W. Strauch, Equipartitions and a distribution
  for numbers: A statistical model for {Benford's} law, Phys. Rev. E 91 (2015)
  062138.
\newblock \href {https://doi.org/10.1103/PhysRevE.91.062138}
  {\path{doi:10.1103/PhysRevE.91.062138}}.

\bibitem{hill1995}
T.~P. Hill, A statistical derivation of the significant-digit law, Statistical
  Science 10~(4) (1995) 354--363.
\newblock \href {https://doi.org/10.1214/ss/1177009869}
  {\path{doi:10.1214/ss/1177009869}}.

\bibitem{ENGEL2003}
H.-A. Engel, C.~Leuenberger, Benford's law for exponential random variables,
  Statistics \& Probability Letters 63~(4) (2003) 361--365.
\newblock \href {https://doi.org/10.1016/S0167-7152(03)00101-9}
  {\path{doi:10.1016/S0167-7152(03)00101-9}}.

\bibitem{proof}
M.~Cong, B.-Q. Ma, A proof of first digit law from {Laplace} transform, Chinese
  Physics Letters 36~(7) (2019) 070201.
\newblock \href {https://doi.org/10.1088/0256-307X/36/7/070201}
  {\path{doi:10.1088/0256-307X/36/7/070201}}.

\bibitem{proof2}
M.~Cong, C.~Li, B.-Q. Ma, First digit law from {Laplace} transform, Physics
  Letters A 383~(16) (2019) 1836--1844.
\newblock \href {https://doi.org/10.1016/j.physleta.2019.03.017}
  {\path{doi:10.1016/j.physleta.2019.03.017}}.

\bibitem{Cigler1964}
J.~Cigler, Methods of summability and uniform distribution mod 1, Compositio
  Mathematica 16 (1964) 44--45,
  \url{http://www.numdam.org/item?id=CM_1964__16__44_0}.

\bibitem{Flehinger1966}
B.~J. Flehinger, On the probability that a random integer has initial digit
  {$A$}, Amer. Math. Monthly 73 (1966) 1056--1061.
\newblock \href {https://doi.org/10.2307/2314636} {\path{doi:10.2307/2314636}}.

\end{thebibliography}






\begin{figure}[H]
    \centering
    \label{graphicalabstract}
   \includegraphics[width=0.8\textwidth]{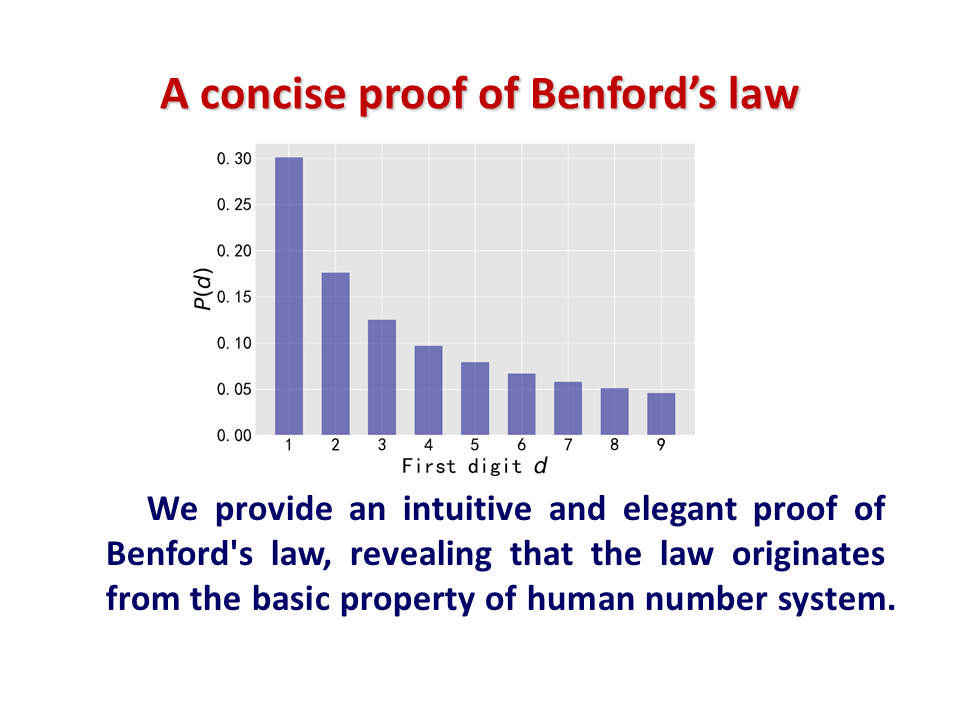} 
    \caption{The graphical abstract of this article, published in Fundamental Research, 4 (2024) 841-844, by L.~Wang \& B.-Q.~Ma,  \href{https://doi.org/10.1016/j.fmre.2023.01.002}{doi:10.1016/j.fmre.2023.01.002}.}
\end{figure}


\end{document}